\documentclass[12pt]{amsart}
\usepackage{amscd,amssymb}
\usepackage[arrow,matrix]{xy}
\usepackage[plainpages,backref,urlcolor=blue]{hyperref}

\theoremstyle{plain}
\newtheorem{thm}[subsection]{Theorem}
\newtheorem{lem}[subsection]{Lemma}
\newtheorem{prop}[subsection]{Proposition}
\newtheorem{cor}[subsection]{Corollary}

\theoremstyle{definition}
\newtheorem{rk}[subsection]{Remark}

\newtheorem{ex}[subsection]{Example}

\numberwithin{equation}{section}
\newcommand{\G}{{\mathcal G  }}
\newcommand{\F}{{\mathbb F}}
\newcommand{\A}{{\mathcal A}}

\newcommand{\al}{{\alpha}}
\newcommand{\be}{{\beta}}

\newcommand{\Z}{\mathbb{Z}}
\newcommand{\Q}{\mathbb{Q}}
\newcommand{\R}{\mathbb{R}}
\newcommand{\C}{\mathbb{C}}

\newcommand{\PP}{\mathbb{P}}
\newcommand{\T}{\mathbb{T}}
\newcommand{\N}{\mathbb{N}}

\DeclareMathOperator{\codim}{codim}

\begin{document}

\title [ Tate properties, polynomial-count varieties, and arrangements ]
{Tate properties, polynomial-count varieties, and monodromy of hyperplane arrangements }

\author[Alexandru Dimca]{Alexandru Dimca$^1$}
\address{  Laboratoire J.A. Dieudonn\'e, UMR du CNRS 6621,
                 Universit\'e de Nice Sophia Antipolis,
                 Parc Valrose,
                 06108 Nice Cedex 02,
                 France}
\email{dimca@unice.fr}

\thanks{$^1$ Partially supported by the  ANR-08-BLAN-0317-02 (SEDIGA)} 



\begin{abstract}

The order of the Milnor fiber monodromy operator of a central hyperplane arrangement is shown to be combinatorially determined.
In particular, a necessary and sufficient condition for the triviality of this monodromy operator is given. 

It is known that the complement of a complex hyperplane arrangement is cohomologically Tate and, if the arrangement is defined over $\Q$, has polynomial count. We show that these properties hold for the corresponding Milnor fibers if the monodromy is trivial.

We construct a  hyperplane arrangement defined over $\Q$, whose Milnor fiber has a nontrivial monodromy operator, is cohomologically Tate, and has not polynomial count. Such examples are shown not to exist in low dimensions.

\end{abstract}

\maketitle

\tableofcontents

\section{Introduction} \label{sec1}

Let $\A$ be a central arrangement of $d$ hyperplanes in $\C^{n+1}$, with $d \geq 2$ and $n\geq 1$, given by a reduced equation 
$ Q(x)=0$. 
Consider the corresponding global Milnor fiber $F$ defined by $Q(x)-1=0$ in $\C^{n+1}$ with monodromy action $h:F \to F$, $h(x)=\exp(2\pi i/d)\cdot x$.

A general investigation line is to check whether certain properties of the associated projective hyperplane arrangement complements $M(\A)$ in $\PP^n$ extend to the Milnor fiber $F$. For instance, note the following.

(i) The cohomology ring $H^*(M(\A),\Z)$ is determined by the combinatorics, see \cite{OS}, but the same question for
the Betti numbers of the Milnor fiber $F$ is widely open, see for instance \cite{Li2}.

(ii) $H^*(M(\A),\Z)$ is torsion free, see \cite{OS}, and there is the open question about the torsion freeness of $H^*(F,\Z)$, see \cite{CDS}. 

(iii) The complement $M(\A)$ is formal in the sense of Sullivan \cite{Sul}, \cite{DGMS}, and it was recently shown that the Milnor fiber $F$ may be not even 1-formal, see \cite{Z}.

\medskip

We say that a complex variety $Y$ is  {\it cohomologically Tate } if for any cohomology group $H^m(Y,\C)$, one has the following vanishing of mixed Hodge numbers: $h^{p,q}(H^m(F,\C))=0$ for $p \ne q$. The fact that the hyperplane arrangement complements $M(\A)$ are cohomologically Tate is known for a long time:  any cohomology group $H^m(M(\A),\Q)$ is a pure Hodge structure of type $(m,m)$, see \cite{L1}, \cite{K}, \cite{Sh}.

  When the monodromy action $h^*$ is {\it trivial on all the cohomology groups} $H^*(F,\C)$, it follows that  we have an equality $H^m(F,\Q)=H^m(M(\A),\Q)$ for any $0 \leq m \leq n$, and hence in this case $F$ is cohomologically Tate. One may ask whether this is the only possibility for a hyperplane arrangement Milnor fiber $F$ to be cohomologically Tate. The claim that $F$ cohomologically Tate implies $h^*$ trivial is shown to be true in the case $n=1$ (obvious, using the MHS on the Milnor fiber of an isolated homogeneous hypersurface singularity, given by Steenbrink in \cite{St1} and recalled in \cite{D0}, pp.243-244) and $n=2$, i.e. for plane arrangements, 
 and negative in general. We do not know whether there is a similar result to Theorem \ref{T1} for $2<n<7$, and this explains why we go to dimension 7 to construct our example.

\medskip

To state our results in this direction, we need some preliminaries.

In studying the cohomology $H^*(F,\Q)$ of the Milnor fiber, the monodromy action 
$h^*: H^*(F,\Q) \to H^*(F,\Q)$ or the number of points in $|F(\F_p)|$, we can, without any loss of generality,  suppose that the arrangement $\A$ is essential, i.e. $\cap_{H\in \A}H=0.$ This is the same as supposing that the polynomial $Q$ involve in an essential way all the variables $x_0,...,x_n$, i.e. one can not choose the coordinates $x$ on $\C^{n+1}$ such that $Q(x_0,...,x_n)=R(x_0,...,x_u)$
for some $0 \leq u <n$ and a homogeneous polynomial $R \in \C[x_0,...,x_u]$ .  

The properties of the monodromy $h^*: H^*(F,\Q) \to H^*(F,\Q)$ are rather misterious, and many things that we know in general are related to the spectrum 
\begin{equation} 
\label{eq2}
Sp(\A)=\sum_{\al \in \Q}m_{\al}t^{\al},
\end{equation}
with $m_{\al}=\sum_j(-1)^{j-n}\dim Gr_F^p\tilde H^j(F,\C)_{\be}$ where $p=[n+1-\al]$ and
$\be=\exp(-2\pi i\al)$, 
which is combinatorially determined, see \cite{BS}.
Surprinsingly, note that for most arrangements the situation is rather simple, namely $h^m: H^m(F,\Q) \to H^m(F,\Q)$ is trivial (i.e. the identity) for $0 \leq m <n$
and $\dim H^n(F,\C)_{\be}=|\chi(M(\A))|$, for any $\be \in \mu_d=\{z \in \C~~|~~z^d=1\}$ with $\be \ne 1$, see for instance \cite{CS}, \cite{Li}, as well as Prop. 2.5.4, Prop.6.4.6, Example 6.4.14 and Theorem 6.4.18 in \cite{D2}. 

\medskip

\begin{thm} \label{T1} 
Let $\A$ be an essential central arrangement of $d$ planes in $\C^3$. The following conditions are equivalent.

\noindent(i) The mixed Hodge numbers  $h^{p,q}(H^2(F,\C))$ vanish for $p \ne q$.

\noindent(ii) The arrangement $\A$ is reducible. 

\noindent(iii) The monodromy action $h^*$ is trivial on all the cohomology groups $H^*(F,\C)$.

\noindent(iv) The spectral numbers $m_{\al}$ in $Sp(\A)$ vanish for all $\al \in (0,1)$.

\end{thm}

Geometrically, the property $(ii)$  means the following: in the projective line arrangement $\A'$ associated to $\A$,
 $(d-1)$ lines meet in one point, say $A$, and the remaining line
$L_d$ does not contain $A$.

In terms of coordinates, this means that 
one may choose the coordinates $(x:y:z)$ on $\PP^2$ such that $A=(0:0:1)$ and $L_d:z=0$. With this choice one has $Q(x,y,z)=Q_1(x,y)z$, where $Q_1$
is a degree $(d-1)$ reduced homogeneous polynomial in $x,y$. This property is exactly the definition of a reducible arrangement when $n=2$.

\medskip

To discuss the higher dimensional case we need the following precise characterization of arrangements with a trivial monodromy.

\begin{thm} \label{Thm1} 
For an essential central arrangement $\A$, the following conditions are equivalent.

\noindent(i) The monodromy action $h^*$ is trivial on all the cohomology groups $H^*(F,\C)$.

\noindent(ii) The arrangement $\A$ is reducible and satisfies the following: if $\A=\A_1 \times ...\times \A_q$ is the decomposition of $\A$ as a product of irreducible arrangements and if $d_j=|\A_j|$ denotes the number of hyperplanes in $\A_j$, then $G.C.D.(d_1,...,d_q)=1.$

Moreover, if $\A$ is defined over $\Q$ (i.e. each hyperplane in $\A$ is defined over $\Q$), then all the irreducible arrangements $\A_j$ are also defined over $\Q$.

\end{thm}
We show below, see Lemma \ref{comb}, that in the (unique) decomposition  $\A=\A_1 \times ...\times \A_q$  of $\A$ as a product of irreducible arrangements $\A_j$, the integers $q$, $d_1$,...,$d_q$ are determined by the combinatorics, i.e. by the intersection lattice $L(\A)$, see \cite{OT}. 

We recall that a central arrangement $\A$ as above is {\it reducible} if one can choose the coordinates $x$ on $\C^{n+1}$ such that $Q(x_0,...,x_n)=R_1(x_0,...,x_u)R_2(x_{u+1},...,x_n)$
for some $0 \leq u <n$ and homogeneous non-constant polynomials $R_1$ and $R_2$. We write then: $\A=\A_1 \times
\A_2$, with $\A_j:R_j=0$.

It is known that an essential arrangement $\A$ is reducible if and only if $\chi(M(\A))=0$, see \cite{STV}.
On the other hand, a trivial monodromy action $h^*$ implies $\chi(M(\A))=0$, as a simple consequence of the following general formula 
\begin{equation} 
\label{zeta}
\sum_j(-1)^j \dim H^j(F,\C)_{\be}=\chi(M(\A))
\end{equation}
for any $\be \in \mu_d$, see Prop. 2.5.4 and  Prop.6.4.6 in \cite{D2}. See also (1.4.2) in \cite{BS}.

\begin{rk}
\label{rkForm} All compact K\"ahler manifolds are formal spaces, see \cite{DGMS}.
When the monodromy $h^*$ is trivial, the corresponding Milnor fiber is clearly a formal space in Sullivan's sense, see
\cite{Sul}, \cite{DP}. Hence Theorem \ref{Thm1} yields a wealth of new examples of formal spaces in the class of smooth affine varieties. Note also that there are examples of non-formal smooth affine surfaces, either related to Milnor fibers of central plane arrangements, see \cite{Z}, or to isolated weighted homogeneous singularities \cite{DPS}.

\end{rk}

Theorem \ref{Thm1} is an obvious consequence of the following {\it Thom-Sebastiani type result},
where the sum of polynomials in disjoint sets of variables is replaced by their product.
This result plays also a key role in the construction of our examples below.

\begin{thm} 
\label{T2}

Let $\A=\A_1 \times ...\times \A_q$ be the decomposition of the central essential arrangement $\A$ as a product of irreducible arrangements, let $d_j=|\A_j|$ denotes the number of hyperplanes in $\A_j$ and let $d_0=G.C.D.(d_1,...,d_q)$. Then the following hold.

\noindent(i) There is a natural identification of graded MHS defined over $\R$
$$H^*(F,\C)=H^*(\T,\C) \otimes (\oplus_{\be \in \mu_{d_0}}(H^*(F_1,\C)_{\be} \otimes \cdots  \otimes H^*(F_q,\C)_{\be})).$$
More precisely, for any $\be \in \mu_{d_0}$, there is an identification
$$H^*(F,\C)_{\be}=H^*(\T,\C) \otimes H^*(F_1,\C)_{\be} \otimes \cdots  \otimes H^*(F_q,\C)_{\be}.$$

\noindent(ii) The monodromy operator $h^*: H^*(F,\C) \to H^*(F,\C)$ has order ${d_0}$, which is determined by the intersection lattice $L(\A)$.

\end{thm} 

Since each $\A_j$ is irreducible, it follows from \eqref{zeta} that each $H^*(F_j)_{\be}$ which occur in Theorem \ref{T2} is nonzero. 

\begin{rk} 
\label{rkT}
There are a number of papers dealing with Thom-Sebastiani type results for the product of two polynomials $f$ and $g$ (or, more generally, for $h(f,g)$, with $h$ a function of two variables), see for instance Oka \cite{O},
Sakamoto \cite{Sk}, N\'emethi \cite{N} and Tapp \cite{T}. The decomposition in our Theorem \ref{T2}, (i), appears in Tapp \cite{T}. However, in all of these papers there is no reference to the mixed Hodge structures
involved, and the methods used, being purely topological, do not allow in a direct way to derive such conclusions. Since these MHS play a central role in the sequel of our paper, we give below a proof
of this decomposition rather different from that in \cite{T}, carefully handling the corresponding MHS.

\end{rk} 

\medskip

We have noticed already that $F$ is a cohomologically Tate variety when the monodromy action $h^*$ is trivial. We give below an example showing that the converse claim is false in general, see Example \ref{ex1}.
This Example is also interesting since it shows that the part $H^{<top}(F,\C)_{\ne 1}$ of the Milnor fiber cohomology can be rather big, unlike all the previously known examples.

\begin{cor}
\label{cor2} Consider the central rational hyperplane arrangement $\A_{u,v}$ in $\C^{n+1}$ defined in Example \ref{ex1},
for any $u,v \in \Z_{>0}$ . Then $n=3u+5v-1$, the only eigenvalues of the monodromy operator
$h^*$ are $\pm 1$ and $\dim H^*(F_{u,v},\C)_{-1} =2^{u+v-1}$. More precisely, the nontrivial $(-1)$-eigenspaces are exactly $H^{2u+4v+j}(F_{u,v},\C)_{-1}$ for $0 \leq j \leq u+v-1$
and
$$\dim H^{2u+4v+j}(F_{u,v},\C)_{-1}={u+v-1 \choose j}.$$

\end{cor}

\medskip

The smallest $n$ for which our construction yields a counterexemple is $n=7$; the Milnor fiber of the corresponding hyperplane arrangement  $\A_{1,1}$ is cohomologically Tate, but $h^*$ is not trivial.

\bigskip

\noindent Recall  that the (compactly supported) Hodge-Deligne polynomial (or E-polynomial) associated to a complex variety $Y$ is given by
\begin{equation} 
\label{HD1}
HD_X(x,y)=\sum_{u,v}(\sum_j(-1)^jh^{u,v}(H^j_c(Y,\C)))x^uy^v. 
\end{equation}
Assume that $Y$ is in fact defined over $\Q$. We say that $Y$  has {\it polynomial count with count polynomial} $P_Y$, if there is a polynomial $P_Y\in \Z[t]$, such that for all but finitely many primes $p$, for any finite field
$\F_q$ with $q=p^s$ and  $s \in \N^*$, the number of points of $Y$ over $\F_q$ is precisely $P_Y(q)$.
As an example, if the hyperplane arrangement $\A$ is
{\it  defined over $\Q$ }, then $M(\A)$  has polynomial count, see for instance \cite{DL0}, section (5.3) or \cite{Stan}, Theorem 5.15. The above rationality assumption is essential: indeed, the line arrangement $\A: x^2+y^2=0$ is defined over $\Q(i)$, and the corresponding complement $M(\A)$ satisfies $|M(\A)(\F_p)|=p-1$ if $p \equiv 1$ mod $4$ and
$|M(\A)(\F_p)|=p+1$ if $p \equiv 3$ mod $4$. Note that in this case $Q(x)$ has integer coefficients, $\A$ is reducible, but the corresponding splitting $\A =\A_1 \times \A_2$ is not defined over $\Q$.

\medskip

A theorem of N. Katz  in \cite{HR} says that if $Y$ has polynomial count with count polynomial $P_Y$, then $HD_Y(x,y)=P_Y(xy)$. See Theorem 2.1.8 in \cite{HR}, and the remark at the bottom of page 563 and Example 2.1.10 explaining that there can be allowed finitely many 'bad characteristics' $p$.
In particular, such a variety is not too far from being cohomologically Tate, i.e. the not Tate part of the cohomology should cancel out in $HD_Y(x,y)$.

\medskip

One may ask what happens to the Milnor fibers of hyperplane arrangements. Again we need a rationality assumption:
the Milnor fiber $F:x^2+y^2=1$ satisfies $|F(\F_p)|=p-1$ if $p \equiv 1$ mod $4$ and
$|F(\F_p)|=p+1$ if $p \equiv 3$ mod $4$. Note that in this example the monodromy operator $h^*$ is trivial and
$F$ is cohomologically Tate.

When $\A$ is defined over $\Q$, we will always choose the defining equation $Q(x)=0$ with integer coefficients.
It turns out that the arithmetic properties of the corresponding Milnor fiber $F:Q(x)-1=0$ do not depend on the choice of $Q$, see for instance Corollary \ref{corT}.

So a first naive idea when trying to construct cohomologically Tate Milnor fibers of rational arrangements which have not polynomial count is to look for a hyperplane arrangement $\A$ with a trivial monodromy $h^*$ and such that the corresponding Milnor fiber $F$ has not polynomial count. Such an attempt cannot succeed in view of the following result, implied by Theorem  \ref{Thm1}.

\begin{cor}
\label{corT} 
If the monodromy action $h^*$ is trivial on all the cohomology groups $H^*(F,\C)$
of the Milnor fiber $F$ of an essential central arrangement $\A$ defined over $\Q$, then $F$ has polynomial count, with the same count polynomial as the corresponding projective complement $M(\A)$.
\end{cor}
In particular, this property of the Milnor fiber of having polynomial count, and the corresponding count polynomial, do not depend on the choice of defining equation $Q$ in $\Z[x]$.

For $n=2$ we have the following.

\begin{cor}

\label{corT2} Consider the following conditions.

\noindent(i)  $Y$ has polynomial count.

\noindent(ii) $Y$ is cohomologically Tate.

If $Y$ is a smooth affine surface, then $(i) \Longrightarrow (ii)$.
If in addition $Y$ is the Milnor fiber of a central rational plane arrangement in $\C^3$, then one also has $(ii) \Longrightarrow (i)$.

\end{cor}

The hyperplane arrangement  $\A_{1,1}$ introduced in Example \ref{ex1} is defined over $\Q$, and the corresponding  Milnor fiber is cohomologically Tate, but has not polynomial count, see Theorem \ref{T5}.

\bigskip

We would like to thank Nero Budur, Denis Ibadula, Mark Kisin, Gus Lehrer and Morihiko Saito for useful discussions in relation to this
work. Special thanks are due to Gabriel Sticlaru who provided expert help in some numerical experiments
with rather large numbers, see Remark \ref{rk2}.

\section{Proof of Theorem \ref{T2} and of Corollary \ref{corT}} \label{sec2}

\proof

We can always, up to a linear coordinate change, write $Q(x)=Q_1(y_1)\cdots Q_q(y_q)$, where $x=(y_1,...,y_q) \in \C^{n+1}$, $y_j \in \C^{n_j}$
such that: $n_1+ ...+n_q=n+1$, $d_j=\deg Q_j>0$ and $\A_j:Q_j=0$ an irreducible and essential arrangement in $\C^{n_j}$ for $j=1,...,q.$

Note that the existence of such a decomposition is equivalent to the following property: there is a partition
of the hyperplanes in $\A$ in $q$ subsets $A_1$,...,$A_q$, such that if we define $V_j=\cap_{H \notin A_j}H$,
then there is a direct sum decomposition $\C^{n+1}=V_1+V_2+...+V_q$, and this partition is the finest with this property. More precisely we have the following result.

\begin{lem} 
\label{comb} Let $\A=\{H_i\}_{i \in I}$ be a central, essential hyperplane arrangement in $\C^{n+1}=V$.
Consider the set $P$ of partitions $I=I_1\cup...\cup I_m$ of $I$ satisfying the following condition:
$\sum_j\codim (\cap_{i\in I_j}H_i)=\dim V$.
Then the following hold.

(i) $P$ is non-empty, since the trivial partition $I=I$ is in $P$.

(ii) If $I=I_1\cup...\cup I_m$ and $I=I'_1\cup...\cup I'_{m'}$ are two partitions in $P$, then their intersection $I=\cup_{i,j}I_{i,j}$ where $I_{i,j}=I_i\cap I_j$ for $i=1,m$,  $j=1,m'$ (the empty intersections $I_{i,j}$ are discarded), is again a partition in $P$.

(iii) In the (unique) decomposition  $\A=\A_1 \times ...\times \A_q$  of $\A$ as a product of irreducible arrangements $\A_j$, the integers $q$, $d_1$,...,$d_q$ are determined by the combinatorics, i.e. by the intersection lattice $L(\A)$.

\end{lem}

\proof

The claim $(i)$ is obvious.

Let $E$ be the dual of $V$ and, for a given partition $I=I_1\cup...\cup I_m$ in $P$, let $E_i$ be the vector subspace in $E$ spanned by the equations of the hyperplanes in $I_i$. Since $\A$ is essential, it follows that $\sum_jE_j=E$. Since $\dim E_j=\codim (\cap_{i\in I_j}H_i)$, it follows that the above sum is in fact a direct sum.

When we have two partitions as above, we define $E_{i,j}$ to be the vector subspace in $E$ spanned by the equations of the hyperplanes in $I_i\cap I'_j$. Then it follows immediately that the sum $\sum_jE_{i,j}=E_i$
is again direct sum for all $i$'s. Hence the sum $\sum_{i,j}E_{i,j}=V$ is a direct sum, which proves the claim $(ii)$.

Since the set of partitions $P$ is defined only in terms of the intersection lattice $L(\A)$, it follows that the unique minimal element of $P$ (with respect to the partial order given by refinement), is combinatorially determined. This minimal element is denoted above by $I=A_1 \cup ...\cup A_q$.
The direct sum decomposition $\C^{n+1}=V_1+V_2+...+V_q$ mentionned above is just the dual of the direct sum decomposition $\sum_jE_j=E$. Since $d_j=|A_j|$, they are also determined by the combinatorics.

\endproof

\begin{rk}
\label{rkQ}
When $\A$ is defined over $\Q$, it follows that all the vector subspaces $V_j$ are defined over $\Q$.
Then the arrangement $\A_j$, which is essentially given by the traces of $H \in A_j$ on $V_j$, is clearly defined over $\Q$. Moreover, the coordinate change from the coordinates $x$ to the coordinates $y$ is defined over $\Q$. This shows that when doing computations for almost all primes $p$, we may replace the equation
$Q(x)=0$ (resp. $Q(x)=1$) by the corresponding equations $Q_1(y_1)\cdots Q_q(y_q) =0$ (resp. $Q_1(y_1)\cdots Q_q(y_q) =1$).

\end{rk}

\medskip

{\bf Proof of Theorem \ref{T2}.}

Let $F_j:Q_j=1$ and $h_j:F_j\to F_j$ be the corresponding Milnor fibers and monodromy homeomorphisms. Let us consider the corresponding least common multiple $m=L.C.M.(d_1,...,d_q)$ and set $w_j=m/d_j$ for  $j=1,...,q.$

Our first aim is to obtain a description of the (total) Milnor fiber $F$ in terms of the
collection of Milnor fibers $F_1$,...,$F_q$.
For this we consider the affine torus
\begin{equation} 
\label{key1}
\T=\{t=(t_1,...,t_q) \in (\C^*)^q~~|~~t_1t_2\cdots t_q=1\}.
\end{equation}

Consider the mapping
\begin{equation} 
\label{key2}
 f:\T \times F_1 \times \cdots \times F_q \to F
\end{equation} 
given by 
\begin{equation} 
\label{key3}
(t,y_1,...,y_q) \mapsto (t_1^{w_1}y_1,...,t_q^{w_q}y_q).
\end{equation}
It is easy to check that this mapping $f$ is surjective and one has 
$$f(t,y_1,...,y_q)=f(t',y'_1,...,y'_q)$$
if and only if the points $(t,y_1,...,y_q)$ and $(t',y'_1,...,y'_q)$ are in the same $G$-orbit,
where the group 
\begin{equation} 
\label{key4}
G=\{g=(g_1,...,g_q)\in \mu_m^q~~|~~g_1g_2 \cdots g_q=1\}
\end{equation}
acts on $X= \T \times F_1 \times \cdots \times F_q$ via
\begin{equation} 
\label{key5}
g \cdot ((t_1,...,t_q),y_1,...y_q)=((g_1^{-1}t_1,...,g_q^{-1}t_q),g_1^{w_1}y_1,...,g_q^{w_q}y_q) .
\end{equation}

It follows that $F=X/G$ and in particular $H^*(F,\Q)=H^*(X,\Q)^G$, the $G$-fixed part of the cohomology of $X$ under the induced $G$-action. This is an isomorphism of MHS (mixed Hodge structures),
since the $G$-action is algebraic.
Note that 
\begin{equation} 
\label{key5.5}
H^*(X,\C)=H^*(\T,\C) \otimes H^*(F_1,\C) \otimes \cdots  \otimes  H^*(F_q,\C).
\end{equation} 
Moreover, the group $G$ acts trivially on the factor $H^*(\T,\C)$, since $\T$ is a connected algebraic group and $G \subset \T$.
If we set
$$H^*=H^*(F_1,\R) \otimes \cdots  \otimes  H^*(F_q,\R)$$
then it follows that
$$H^*(F,\R)=H^*(\T,\R) \otimes (H^*)^G.$$
The $G$-action on $H^*_{\C}$, the complexification of $H^*$, is given by the following: 
if $\eta= \eta_1 \otimes \cdots \otimes \eta_q$, then
\begin{equation} 
\label{key8.1}
g\eta=   (h_1)^{k_1}(\eta_1) \otimes \cdots  \otimes (h_q)^{k_q}(\eta_q).
\end{equation}
Here $g=(\lambda^{k_1},...,\lambda^{k_q})$ with $\lambda=\exp(2\pi i/m)$ and $k_1+...+k_q$ is divisible by $m$, the $k_j$ being otherwise arbitrary integers.

Let $\eta_j \in H^*(F_j,\C)$ be now chosen such that for any $j=1,...,q$ there is a $\be_j \in \mu_{d_j} \subset \mu_m$ with $h_j^*\eta_j=\be_j\eta_j$ and look at $\eta= \eta_1 \otimes \cdots \otimes \eta_q$. Such elements form a $\C$-basis of $H^*_{\C}$ and hence to determine $(H^*_{\C})^G$ is the same as finding all
$\eta$'s of this form which are fixed under the $G$-action. By choosing $k_q=m-k_1-...-k_{q-1}$
we get from \eqref{key8.1} the following
$$g\eta=(\frac{\be_1}{\be_q})^{k_1} \cdots (\frac{\be_{q-1}}{\be_q})^{k_{q-1}}\eta$$
where now there is no condition on the integers $k_1$,...,$k_{q-1}$. By taking one of them equal to 1 and the rest zero, we see that $\eta \in (H^*_{\C})^G$ implies $\be_1=...=\be_q$. Call this common value $\lambda_0$ and note that $\lambda_0 \in \mu_{d_0}=\cap _{j=1,q}\mu_{d_j}$.

Conversely, for any  $\lambda_0 \in \mu_{d_0}$ and any 
$\eta \in H^*(F_1,\C)_{\lambda_0} \otimes \cdots  \otimes H^*(F_q,\C)_{\lambda_0}$, we see by using \eqref{key8.1} that $\eta \in (H^*_{\C})^G$.

Moreover, set $E_{\lambda_0}=H^*(F_1,\C)_{\lambda_0} \otimes \cdots  \otimes H^*(F_q,\C)_{\lambda_0}$ and note that ${\overline E_{\lambda_0}}=E_{\overline \lambda_0}$.

It follows that if we set $M_{\lambda_0}=E_{\lambda_0}$ when ${\lambda_0} \in \R$ and
$M_{\lambda_0}=E_{\lambda_0}+ E_{\overline \lambda_0} $ when ${\lambda_0} \notin \R$,
then $M_{\lambda_0}$ is endowed with a natural $\R$-MHS, see also the formula \eqref{mhs1}.

\medskip

To prove the second part of the claim $(i)$ in Theorem \ref{T2}, we construct a nice $G$-equivariant lifting $\tilde h:X \to X$ of the monodromy morphism $h:F\to F$. We set
\begin{equation} 
\label{lift1}
\tilde h(t_1,...,t_q,y_1,...,y_q)=(\gamma_1t_1,...,\gamma_qt_q, \be_1y_1,y_2,...,y_q)
\end{equation}
where $\be_1=\exp(2\pi i/d_{1})$ and $\gamma_j=\exp(2\pi i a_j)$ with $w_1a_1= 1/d -1/d_{1}$ and
$w_jaj=1/d$ for $j>1$. Then $\sum_ja_j=0$, i.e. $(\gamma_1t_1,...,\gamma_qt_q)\in \T$ and $f \circ \tilde h=h\circ f$.
Then $\tilde h ^*$ acts as identity on all the factors in the tensor product \eqref{key5.5}, except
on $H^*(F_1,\C)$, where it acts via $h_1^*$.
We get the second part of the claim $(i)$ by using the description of the cohomology $H^*(F,\C)$ given in the first part of the claim $(i)$.

\bigskip

To prove the claim $(ii)$, note that the affine torus $\T$ acts on the Milnor fiber $F$ by 
\begin{equation} 
\label{key9}
t(y_1,...,y_q)=(t_1^{w_1}y_1,...,t_q^{w_q}y_q).
\end{equation}
Hence to show that $(h^*)^{d_0}$ is trivial, it is enough to show the existence of an element
$t=(t_1,...,t_q) \in \T$ such that $t_j^{w_j}=\exp(2\pi id_0/d)$ for $j=1,...,q$.
Since $G.C.D.(d_1,...,d_q)=d_0$, there are integers $k_j$ such that
\begin{equation} 
\label{key10}
k_1d_1+...k_qd_q=(m-1)d_0.
\end{equation}
For $j=1,...,q$ we set
\begin{equation} 
\label{key11}
t_j=\exp[2\pi i(\frac{d_0}{dw_j}+\frac{k_j}{w_j})].
\end{equation}
The relations $t_j^{w_j}=\exp(2\pi id_0/d)$ are clearly satisfied. Moreover
$$t_1t_2\cdots t_q=\exp[2\pi i(\sum_j(\frac{d_0d_j}{dm}+\frac{k_jd_j}{m}))]=\exp[2\pi id_0(\frac{1}{m}+\frac{m-1}{m})]=1.$$
Hence $t \in \T$, it follows that $(h^*)^{d_0}$ is trivial. We conclude using following fact: since each $\A_j$ is irreducible, it follows from \eqref{zeta} that each $H^*(F_j)_{\be}$ which occur in Theorem \ref{T2} is nonzero. Hence each $H^*(F)_{\be}$ is nonzero, i.e. the order of $h^*$ is indeed $d_0$.

\endproof

\begin{rk}
\label{rk1} Consider the central essential hyperplane arrangement $\A$ in $\C^{n+1}$ and its decomposition $\A=\A_1 \times ....\A_q$ as a product of irreducible arrangements $\A_j$ for $j=1,...,q$. Then it is easy to show that
$$M(\A)=\T \times M(\A_1)\times ...\times M(\A_q).$$
This implies the following for the cohomology of the corresponding projective complements 
$$H^*(M(\A))=H^*(\C^*)^{\otimes (q-1)} \otimes H^*(M(\A_1)) \otimes...\otimes H^*(M(\A_q))$$
i.e. the case $\be=1$ in Theorem \ref{T2}, $(iii)$. See also Tapp \cite{T}.

\end{rk}

Now we pass to the proof of Corollary \ref{corT}. Since the monodromy $h^*$ is trivial, we know by Theorem \ref{Thm1} that $G.C.D.(d_1,...d_q)=1$. Hence there exist integers $m_j$ such that $m_1d_1+...m_qd_q=1$. 

Let $\F$ be a finite field and consider the mapping $Q:\F^{n+1} \to \F$ induced by the polynomial $Q$.

For $a \in \F$, denote by $F(a)$ the fiber $Q^{-1}(a)$. Denote also $M(\A,\F)$
the corresponding (projective) hyperplane arrangement complement over $\F$. It is clear that
\begin{equation} 
\label{t1}
|\F^{n+1} \setminus F(0)|=(|\F|-1)\cdot |M(\A,\F)|.
\end{equation}
Consider the following $\F^*$-action on $\F^{n+1}$:
\begin{equation} 
\label{t2}
t\cdot x= t \cdot(y_1,y_2,...,y_q)=(t^{m_1}y_1,t^{m_2}y_2,...,t^{m_q}y_q)  .
\end{equation}
The relation $Q(t\cdot x)=t Q(x)$ shows that all the fibers $F(a)$ for $a \in \F^*$ have the same cardinal. Since their disjoint union is exactly $\F^{n+1} \setminus F(0)$, the equation
\eqref{t1} yields
\begin{equation} 
\label{t3}
|F(1)|= |M(\A,\F)|.
\end{equation}
This equality completes the proof of Corollary \ref{corT}.

\section{Proof of Theorem \ref{T1} and of Corollary \ref{corT2}} \label{sec3}

\proof

It follows from the discussion just after Theorem \ref{T1} that, assuming $(ii)$, the Milnor fiber $F$ is isomorphic to the complement of the central line arrangement given by $Q_1=0$ in $\C^2$
(indeed, the only partitions of $3$ are $1+2=2+1=1+1+1=3$).
Hence the implication $(ii) \Rightarrow (i)$ is obviously true.

Note that for $n=2$ and $\al \in (0,1)$, the corresponding spectral number is by definition
\begin{equation} 
\label{k1}
m_{\al}=h^{2,0}(H^2(F,\C)_{\beta})+h^{2,1}(H^2(F,\C)_{\beta})
\end{equation}
with $\beta=\exp(-2\pi i \al)$. Indeed, $h^{2,2}(H^2(F,\C)_{\beta}=0$, as follows from Theorem 1.3 in \cite{DL}. 

Since $H^2(F,\C)_1$ is known to be of type $(2,2)$, the equivalence of the claims $(i)$ and $(iv)$ in Theorem \ref{T1} follows. Moreover, the equivalence of the claims $(ii)$ and $(iii)$ follows from
Theorem \ref{Thm1}, using again the partions of $3$.

\bigskip

So from now on we assume that $(iv)$ holds and we prove $(ii)$. The fact that the arrangement is essential
implies that $d \geq 3$ and that the lines in $\A'$ do not pass all through the same point, i.e. there is no point $s$ of multiplicity $m_s=d$. To prove $(ii)$ we have to show the existence of a point of multiplicity $d-1$.

For $d=3$ there are only two type of arrangements, described better in terms of their associated projective line arrangements $\A'$:

(a) three lines forming a triangle, in which case one may take $Q=xyz$ and $F=\C^* \times \C^*$, and

(b)  three lines meeting at one point.

The claim $(i) \Rightarrow (ii)$ is clear by the previous remark.

\bigskip

We assume from now on that $d \geq 4$. 

Next we recall the following key formula from \cite{BS}, Theorem 3 (rewritten slightly for our needs). If $0<\al=\frac{j}{d}<1$, then 
\begin{equation} 
\label{k2}
 m_{\al}={j-1 \choose 2}-\sum_{s; m_s\geq 3}{\lceil jm_s/d \rceil -1 \choose 2},
\end{equation} 
where the sum is over all multiple points $s$ in $\A'$ with multiplicity $m_s \geq 3$. By convention ${a \choose b}=0$ if $a<b$.

If we use the above formula for $j=3$, we get that the corresponding
vanishing $m_{\al}=0$ is equivalent to the existence a unique point $s$ of multiplicity $m_s>2d/3$
in $\A'$. For $d=4$ (resp. $d=5$), this means a point of multiplicity $m_s \geq 3$ (resp. $m_s \geq 4$).
As above (case $d=3$, (b)), the case $m_s=4$ (resp. $m_s=5$) is discarded since $\A'$ is essential. Hence $m_s=3$ (resp. $m_s=4$),
which gives exactly an arrangement $\A'$ as claimed in $(ii)$.

\bigskip

From now on we assume $d > 5$. We apply  the formula \eqref{k2} for $j=d-1$.
Since one clearly has
$$m-1<\frac{(d-1)m}{d}<m,$$
it follows that $\lceil (d-1)m/d \rceil=m$ and hence the vanishing $m_{\al}=0$ in this case is equivalent to the equality
\begin{equation} 
\label{k3}
{d-2 \choose 2}=\sum_{s; m_s\geq 3} {m_s-1 \choose 2}.
\end{equation}
Similarly, for $j=d-2$ we get from $m_{\al}=0$  the following equality
\begin{equation} 
\label{k4}
{d-3 \choose 2}=\sum_{s; 3\leq m_s<d/2} {m_s-1 \choose 2}+ \sum_{s; m_s\geq d/2} {m_s-2 \choose 2}.
\end{equation}
By taking the difference of \eqref{k3} and \eqref{k4} we get the equality
\begin{equation} 
\label{k8}
d-3=\sum_{s; m_s\geq d/2}(m_s-2).
\end{equation}
It follows that, if the set $S=\{s; m_s\geq d/2\}$ contains exactly one element, then the corresponding multiple point $s$
satisfies $m_s=d-1$ and we are done. 

Suppose now that the set $S$ contains at least two elements. Since one of them has to be the point $s$ with multiplicity $m_s>2d/3$ obtained above for $j=3$, we get in this case
\begin{equation} 
\label{k9}
d-3> (2d/3-2)+(d/2-2).
\end{equation}
This is equivalent to $d<6$, a contradiction with our hypothesis $d >5$.

\endproof

The proof of Corollary \ref{corT2} is very easy now. 
If $(i)$ holds, it follows from Katz' Theorem that the Hodge-Deligne polynomial of $Y$ contains only the monomials $1, xy, (xy)^2$. The possible non-zero mixed Hodge numbers in this situations are:
$$h^{2,2}(H^4_c(Y,\C)),~~
h^{1,2}(H^3_c(Y,\C)),~~h^{2,1}(H^3_c(Y,\C)),~~h^{1,1}(H^3_c(Y,\C)),~~
h^{2,0}(H^2_c(Y,\C)),$$
$$h^{1,1}(H^2_c(Y,\C)),~~h^{0,2}(H^2_c(Y,\C)),~~h^{0,1}(H^2_c(Y,\C)),~~h^{1,0}(H^2_c(Y,\C)),~~h^{0,0}(H^2_c(Y,\C)).$$
It follows that the only cancellations in the Hodge polynomial can occur in the coefficient of $xy$.
Hence all the mixed Hodge numbers $h^{u,v}(H^m_c(Y,\C))$ vanish for $u\ne v$, i.e. $Y$ is cohomologically Tate.

Now for Milnor fibers  of a central plane arrangement we have seen in Theorem \ref{T1} that $(ii)$ implies  that  $Y$ is isomorphic to a line arrangement complement in $\C^2$, and hence it has polynomial count.

\section{A purity result and the key  example} \label{sec4}

Let $\A$ be a central arrangement of $d$ hyperplanes in $\C^{n+1}$, with  $n\geq 1$, given by a reduced equation $Q(x)=0$. Then clearly $H^n(F,\Q)_1$ and $H^{n}(F,\C)_{-1}$ are mixed Hodge substructures in $H^n(F,\Q)$. Moreover, for $\be \in \mu_d$, $\be \ne \pm 1$, the same is true for the subspace 
\begin{equation} 
\label{mhs1}
H^n(F,\C)_{\be, \overline \be}=H^n(F,\C)_{\be} \oplus H^n(F,\C)_{ \overline \be}=\ker[(h^n)^2-2Re(\be)h^n+Id]
\end{equation}
which is in fact defined over $\R$ (as the last equality shows). For $\be=-1$, we set
$H^n(F,\C)_{\be, \overline \be}=H^{n}(F,\C)_{-1}$ for uniformity of notation.

Let $D= Q^{-1}(0)=\cup_{H \in \A}H$. For a point $x \in D$, $x \ne 0$, let $L_x=\cap_{H \in \A, x \in H}H$ and denote by $\A_x$ the central hyperplane arrangement induced by $\A$
on a linear subspace $T_x$, passing through $x$ and transversal to $L_x$. We may choose
$\dim T_x=\codim L_x$ and identify $x$ with the origin in the linear space $T_x$.
Let $h^*_x:H^*(F_x, \C) \to H^*(F_x, \C)$ be the corresponding monodromy operator at $x$.

With this notation we have the following result.

\begin{prop}
\label{prop1} 
Let $\be \in \mu_d$, $\be \ne 1$ be a root of unity which is not an eigenvalue for any monodromy operator $h^*_x$ for $x \in D$, $x \ne 0$. Then the corresponding eigenspace $H^n(F,\C)_{\be, \overline \be}$ is a pure Hodge structure of weight $n$. 

In particular, if $\be=\exp(-2\pi i\al)$ for some $\al \in \Q$, then the coefficients in the corresponding  spectrum $Sp(\A)$ have the following symmetry property: 
\begin{equation} 
\label{e2}
m_{\al}=m_{n+1-\al}.
\end{equation}

\end{prop}

\proof
This result is a direct consequence of Lemma 3.6 in \cite{Sa1}.  
Indeed, our hypothesis on $\be$ implies that the nearby cycle sheaf $\psi_{Q,\be}\C$  is supported at the origin and hence it is identified to $H^{n}(F,\C)_{\be}$. This identification in turn implies that
the logarithm of the unipotent part of the monodromy $N$ is trivial on $\psi_{Q,\be}\C$ (as this holds for $H^{n}(F,\C)_{\be}$,  the monodromy $h^*$ being semisimple).

The last claim about the symmetry property in \eqref{e2} is proved in the usual way.  In view of the purity result, we have $h^{p,q}(H^n(F,\C)_{\be})=0$ for $p+q \ne n$. 

For example, assume that $\al=\frac{j}{d}$ with $0<j<d$.
It follows from the definition of the spectrum \eqref{eq2} that $m_{\al}=\dim h^{n,0}(H^n(F,\C)_{\be})$, $m_{n+1-\al}=\dim h^{0,n}(H^n(F,\C)_{\overline \be})$. The claimed equality follows by considering the action of the complex conjugation on the pure Hodge structure  $H^n(F,\C)_{\be, \overline \be}$.

\endproof

\begin{cor}
\label{cor1} Assume that $\A$ is a generic central hyperplane arrangement, i.e. the associated projective divisor in $\PP^{n}$ is a divisor with normal crossing.
Then 
$$H^n(F,\C)_{\ne 1} =\oplus_{\be \ne 1} H^n(F,\C)_{\be} $$ 
is a pure Hodge structure of weight $n$ and for any rational number
$\al \in \Q \setminus \Z$ one has $m_{\al}=m_{n+1-\al}$.

\end{cor}

This symmetry of the coefficients of the spectrum $Sp(\A)$ of a generic central arrangement is alluded to in \cite{Sa1}, see the remarks just after Corollary 1 in the Introduction. 
The result above follows from Proposition \ref{prop1} using the simple fact that in this case all the monodromy operators $h_x^*$ are the identity. For more on the monodromy of generic arrangements, see pp. 209-210 in \cite{OT}, (3.2) in \cite{BS} and section 3 in \cite{CS}.

Now we pass to our example.

\begin{ex}
\label{ex1}

For $n\geq 2$ let $\G_n$ be the central  arrangement in $\C^{n+1}$ given by the equation
$$Q_n(x)=Q_n(x_0,...,x_n)=x_0\cdot x_1 \cdot ... \cdot x_n \cdot (x_0+x_1+...+x_n).$$
Hence the degree $d_n$ of $Q_n$ is $n+2$ and clearly $\G_n$ is a generic, irreducible arrangement.  Assume that $n=2k$ is even and
use Theorem 3.2 in \cite{CS} (or refer to the original paper \cite{OR}) to get that
\begin{equation} 
\label{e1}
H^m(F,\C)=H^m(F,\C)_1
\end{equation}
for $0 \leq m <n$ and $\dim H^{n}(F,\C)_{-1}=1$.

It follows that all the spectral coefficients in $Sp(\G_n)$ corresponding to the monodromy eigenvalue 
$-1$ vanish except for $m_{k+\frac{1}{2}}=1$ (which is the only auto-dual element in the sum with respect to the symmetry given by Corollary \ref{cor1}).
It follows that the eigenspace $H^{n}(F,\C)_{-1}$ is spanned by a cohomology class $\omega_n$ of Hodge type $(k,k)$.

Let us return now to the setting of the proof of Theorem \ref{Thm1} in Section \ref{sec2}.
Let $\A_{u,v}$ be the central hyperplane arrangement obtained by taking the product of $u>0$ copies of the arrangement $\G_2$ and $v>0$ copies of $\G_4$. In follows that $n=3u+5v-1$,
$d=4u+6v$, $q=u+v$, $d_0=G.C.D.(d_1,...,d_q)=2$.

In this case, the cohomology of the corresponding total Milnor fiber $F_{u,v}$ can be described via Theorem \ref{T2} as the following direct sum $H^*(F_{u,v},\C)=H^*(F_{u,v},\C)_1 \oplus H^*(F_{u,v},\C)_{-1}$, where
\begin{equation} 
\label{e3}
H^*(F_{u,v},\C)_1 =H^*(\T,\C)\otimes H^*(M(\G_2),\C)^{\otimes u} \otimes H^*(M(\G_4),\C)^{\otimes v}
\end{equation}
and
\begin{equation} 
\label{e4}
H^*(F_{u,v},\C)_{-1} =H^*(\T,\C)\otimes (\C \omega_2)^{\otimes u} \otimes (\C \omega_4)^{\otimes v}. 
\end{equation}

It follows that $F_{u,v}$ is a cohomologically Tate variety of dimension $n=3u+5v-1$, but the corresponding monodromy action $h^*$ is not trivial. By choosing various values for $u,v$ we can get $n=7$ as a minimal value as well as any integer $n \geq 15$.

\end{ex}

\begin{rk}
\label{rk21} 
If one is interested only in irreducible arrangements, then such examples with large cohomology $H^{<top}(F,\C)$
can be obtained by taking a generic hyperplane section of the arrangement $\A_{u,v}$.
\end{rk}

\section{Finite field computations} \label{sec5}

In this section we use the following notation. Let $F=F_{1,1}$ be the variety defined over $\Z$ by $Q(x)=1$, where 
$$Q(x)=x_1\cdot x_2\cdot x_3 \cdot (x_1+x_2+x_3)\cdot x_4 \cdot x_5 \cdot x_6 \cdot x_7\cdot x_8\cdot (x_4+x_5+x_6+x_7+x_8).$$
For any prime $p$ we denote by $A(p)$ the number of points of $F$ over $\F_p=\Z/p\Z$, i.e. the number of solutions of the equation $Q(x)=1$ in $\F_p^8$. For $a \in \F_p^*$, consider
the varieties
$$F_{1}(a)=\{(x_1,x_2,x_3)\in \F_p^3~~|~~x_1\cdot x_2\cdot x_3 \cdot (x_1+x_2+x_3)=a\}$$
and
$$F_{2}(a)=\{(x_4,x_5,x_6,x_7,x_8)\in \F_p^5~~|~~ a \cdot x_4\cdot x_5\cdot x_6\cdot x_7\cdot x_8 \cdot (x_4+x_5+x_6+x_7+x_8)=1 \}.$$
Let $n_1(a)=|F_1(a)|$, $n_2(a)=|F_2(a)|$ and note that obviously one has
\begin{equation} 
\label{n1}
A(p)=\sum_{a \in \F_p^*}n_1(a)n_2(a). 
\end{equation}
From now on we assume that $p$ is a prime number in the arithmetic progression $12b+11$, where $b \in \N$. 
\begin{lem} 
\label{lem1}
Consider the group morphism $p(k):\F_p^* \to \F_p^*$ given by $t \mapsto t^k$.
Then the following hold.

(i) There is an index 2 subgroup $H \subset \F_p^*$ such that $p(4)(\F_p^*)=p(6)(\F_p^*)=H$.

(ii) If $a$ and $a'$ have the same class in $\F_p^*/H$, then  $n_1(a)=n_1(a')$ and $n_2(a)=n_2(a')$.
\end{lem}
\proof
Since the multiplicative group $\F_p^*$ is cyclic, of order $p-1=12b+10$, it follows that
to prove (i)
it is enough to show that $p(4)(\F_p^*)$ and $p(6)(\F_p^*)$ have both cardinal $(p-1)/2$.
This is the same as proving that the corresponding kernels have order $2$. Note that the equation $t^4=1$ (resp. $t^6=1$) are both equivalent to $t^2=1$ (look at the order of $t$, which must be a divisor of $12b+10$). Hence $\ker p(4)= \ker p(6)= \{\pm 1\}$ and hence have order $2$.

\bigskip

To prove (ii), consider the action of $\F_p^*$ on $\F_p^3$ (resp. $\F_p^5$) given by the usual multiplication of a vector by a scalar.  This multiplication by $t \in \F_p^*$ induces a bijection $F_1(a)=F_1(t^4a)$ (resp. $F_2(a)=F_2(t^6a)$). This completes the proof.

\endproof

We denote $n_1'=n_1(1)$, $n_2'=n_2(1)$, $n_1''=n_1(a)$ and $n_2''=n_2(a)$, for $a \in (\F_p^* \setminus H)$. We also write $p=12b+11=4k+3$, i.e. we set $k=3b+2$.
With this notation, \eqref{n1} may be rewritten as
\begin{equation} 
\label{n2}
A(p)=|H|(n_1'n_2'+n_1''n_2'')=(2k+1)(n_1'n_2'+n_1''n_2''). 
\end{equation}

The (affine or projective) hyperplane arrangement complements have polynomial count, see for instance \cite{DL0}, section (5.3) or \cite{Stan}, Theorem 5.15. So we apply this fact to the
{\it projective} hyperplane arrangement complements $\A_1'$ in $\PP^2$ (resp. $\A_2'$ in $\PP^4$)
corresponding to $x_1x_2x_3(x_1+x_2+x_3)=0$ (resp. $x_4x_5x_6x_7x_8(x_4+x_5+x_6+x_7+x_8)=0$).
Using the formulas for the Betti numbers in \cite{CS}, the fact that each cohomology group $H^m$ is pure of type $(m,m)$ and the duality between $H^m$ and $H^{2d-m}_c$
(where $d$ is the corresponding (complex) dimension), it follows that

\begin{equation} 
\label{HD2}
HD_{M(\A_1')}(x,y)=(xy)^2-3(xy)+3 
\end{equation}
and
\begin{equation} 
\label{HD3}
  HD_{M(\A_2')}(x,y)=(xy)^4-5(xy)^3+10(xy)^2-10(xy)+5.
\end{equation}
It follows that the corresponding counting polynomials are
\begin{equation} 
\label{C1}
P_{M(\A_1')}(t)=t^2-3t+3 
\end{equation}
and
\begin{equation} 
\label{C2}
P_{M(\A_2')}(t)=t^4-5t^3+10t^2-10t+5.
\end{equation}
Using these two polynomials, it follows that, for almost all $p=4k+3$, one has
\begin{equation} 
\label{n3}
N_1'(p)=|M(\A_1')(\F_p)|=p^2-3p+3 \equiv 4k+3 \text{ mod } 8
\end{equation}
and
\begin{equation} 
\label{n4}
N_2'(p)=|M(\A_2')(\F_p)|=p^4-5p^3+10p^2-10p+5 \equiv 4k+3 \text{ mod } 8.
\end{equation}
Let $\A_1$ (resp. $\A_2$) denote the corresponding central hyperplane arrangements in $\C^3$ (resp. $\C^5$) and
$M(\A_1)$ (resp. $M(\A_2)$) the associated {\it affine} complements.
Then, using \eqref{n3} and \eqref{n4}, we get
\begin{equation} 
\label{n5}
N_1(p)=|M(\A_1)(\F_p)|=(p-1)(p^2-3p+3) \equiv 4k+6 \text{ mod } 8
\end{equation}
and
\begin{equation} 
\label{n6}
N_2(p)=|M(\A_2)(\F_p)|=(p-1)(p^4-5p^3+10p^2-10p+5) \equiv 4k+6 \text{ mod } 8.
\end{equation}
Now, note that $M(\A_1)$ (resp. $M(\A_2)$) is the disjoint union of all the hypersurfaces $F_1(a)$ (resp. $F_2(a)$) for $a \in \F_p^*$. Lemma \ref{lem1}, \eqref{n5} and \eqref{n6}
yield
\begin{equation} 
\label{n7}
n_1'+n_1''=2(p^2-3p+3) \equiv 6 \text{ mod } 8
\end{equation}
and
\begin{equation} 
\label{n8}
n_2'+n_2''=2(p^4-5p^3+10p^2-10p+5) \equiv 6 \text{ mod } 8.
\end{equation}
It follows that $n_1'' \equiv 6-n_1'$ mod $8$ and $n_2'' \equiv 6-n_2'$ mod $8$.
Using \eqref{n2}, we get
\begin{equation} 
\label{n9}
A(p)=(2k+1)(n_1'n_2'+n_1''n_2'') \equiv 2(2k+1)(2+n_1'n_2'-3n_1'-3n_2') \text{ mod } 8. 
\end{equation}

Now we can state and prove our main result of this section.

\begin{thm} \label{T5} 
The variety $F$ satisfies the following.

\noindent(i) $F$ is a cohomologically Tate variety.

\noindent(ii) $F$ has not polynomial count.

\end{thm}

\proof

The proof of the first claim was given already in Example \ref{ex1}. Moreover, it follows from the description of the cohomology given there, that the corresponding Hodge-Deligne polynomial is
\begin{equation} 
\label{HD4}
  HD_{F}(x,y)=(xy)^7-9(xy)^6+36(xy)^5-82(xy)^4+ 119(xy)^3-110(xy)^2+60(xy)-15.
\end{equation}
Assume that (ii) fails, i.e. that $F$ has polynomial count with polynomial $P_F$. Then according to Katz Theorem \cite{HR}, we must have
\begin{equation} 
\label{C3}
 P_{F}(t)=t^7-9t^6+36t^5-82t^4+ 119t^3-110t^2+60t-15.
\end{equation}
We will reach a contradiction by showing that $P_F(p) \ne A(p)$, for infinitely many primes,
namely for all the primes $p$ in the progression $12b+11$ (Dirichlet prime number theorem). Let $p$ be such a prime and write
$p=4k+3$ as above. Then a simple computation (e.g. using Maple) shows that
\begin{equation} 
\label{M1}
P_F(4k+3)=8(2k+1)(1024k^6+2560k^5+2752k^4+1632k^3+584k^2+129k+15).
\end{equation}
Hence, to complete the proof, it is enough to show that $A(p) \not \equiv 0 \text{ mod } 8$, for any prime $p$ as above.
Using the formula \eqref{n9}, this is equivalent to showing
\begin{equation} 
\label{n10}
2+n_1'n_2'-3n_1'-3n_2' \not \equiv 0 \text{ mod } 4. 
\end{equation}
This in turn follows from the following.

\begin{lem} 
\label{lem2}
Let $p$ be a prime number of the form $4k+3$. Then both $n_1'$ and $n_2'$ are divisible by $4$.

\end{lem}

\proof
 Consider the group $G=\{\pm 1\} \times \{1, \tau\}$ acting on the (finite) Milnor fiber $F_1(1)$ by
$$-1\cdot (x_1,x_2,x_3)=(-x_1,-x_2,-x_3)$$
$$\tau \cdot (x_1,x_2,x_3)=(x_2,x_1,x_3)$$
and
$$-\tau \cdot (x_1,x_2,x_3)=(-x_2,-x_1,-x_3).$$
There are two types of $G$-orbits. First we have the orbits $Gx$, corresponding to points $x$
such that $x_1 \ne x_2$. Then the isotropy group $G_x$ is trivial, and the orbit $Gx$ consists of $4$ points. Since we are interested in a computation modulo $4$, we can forget about such orbits.

\bigskip

The second type of orbit corresponds to points $x$ such that $x_1=x_2=t$. Then $x \in F_1(1)$
is equivalent to
\begin{equation} 
\label{E1}
(x_3+t)^2=\Delta (t)
\end{equation}
where $\Delta(t)=t^2+t^{-2}.$ There are two possibilities:

(a) $\Delta(t) \notin H$, i.e. $\Delta(t)$ is not a square in $\F_p^*$. Then the equation \eqref{E1}
is impossible and we do not get a point in $F_1(1)$.

(b) $\Delta(t) $ is  a square in $\F_p^*$. Then the equation \eqref{E1}
has two solutions, namely $x_3=-t+y$ and $x_3=-t-y$, with $y$ satisfying $y^2=\Delta(t) $. In this way we  get two points in $F_1(1)$. Note that $\Delta(t) =\Delta(-t) $, hence for the point $(-t)$ we are exactly in the same situation as for the point $t$. It follows that for each pair $\{t,-t\}$ (i.e. orbit of the obvious $\{\pm 1\}$-action on $\F_p^*$), we get either
0 points in the (finite) Milnor fiber $F_1(1)$, or we get 4 points.
This proves the claim for $n_1'$.

\bigskip

Consider next the group $G=\{\pm 1\} \times \{1, \tau\} \times\{1, \sigma \}$ acting on the (finite) Milnor fiber $F_2(1)$ by
$$-1\cdot (x_4,x_5,x_6,x_7,x_8)=(-x_4,-x_5,-x_6,-x_7,-x_8)$$
$$\tau \cdot (x_4,x_5,x_6,x_7,x_8)  = (x_5,x_4,x_6,x_7,x_8)$$
and
$$\sigma \cdot (x_4,x_5,x_6,x_7,x_8)=(x_4,x_5,x_7,x_6,x_8)$$
(and their obvious consequences).
 Since we are interested in a computation modulo $4$, we can forget about all orbits
corresponding to points having an isotropy group of order at most 2. Let now $x$ be a point such that $G_x=\{1, \tau\} \times\{1, \sigma \}$. Then $x=(s,s,t,t,x_8)$ and $x \in F_2(1)$
is equivalent to
\begin{equation} 
\label{E2}
(x_8+s+t)^2=\Delta (s,t)
\end{equation}
where $\Delta(t)=(s+t)^2+(st)^{-2}.$ Note that $\Delta(s,t) =\Delta(-s,-t) $, hence for the point $(-s,-t)$ we are exactly in the same situation as for the point $(s,t)$. It follows that for each  orbit of the obvious $\{\pm 1\}$-action on $(\F_p^*)^2$, we get either
0 points in the (finite) Milnor fiber $F_2(1)$, or we get 4 points.
This proves the claim for $n_2'$.

\endproof

\begin{rk}
\label{rk2} Surprisingly, it seems that $P_F(p)=A(p)$ for many primes of the form
$p=4k+1$. Indeed, by a computer aided computation, we have obtained the following.
$$P_F(5)=A(5)=11160, ~~  ~~ P_F(13)=A(13)= 30575400~~ ~~,P_F(17)=A(17)=237920544,$$
$$P_F(29)=A(29)=12579682248,~~ ~~P_F(37)=A(37)=74188920024,$$
$$P_F(41)=A(41)=155950465680,~~ ~~ P_F(53)=A(53)=989657318520,$$
$$P_F( 61) = A(61)=  2708356179720,~~ ~~
P_F( 73) = A(73)=    9757738115280.$$
However one has
$$P_F( 89) =    39954467578608 \ne A(89)= 39843984220188,$$
and
$$  P_F( 97) =    73603528860864  \ne A(97)=    72706366451444.$$          

\end{rk}

\begin{rk}
\label{rk0}
M. Kisin and G. Lehrer  have considered a notion of {\it mixed Tate variety} (imposing conditions on the eigenvalues of Frobenius action on $p$-adic \'etale cohomology), see Definition (2.6) in \cite{KL}
and shown that if a variety $Y$ is mixed Tate then  $Y$ is cohomologically Tate (see Theorem (2.2) (2) in \cite{KL}). Then in Remark (2.4) and in a footnote on p. 213, they discuss {\it the conjectural equivalence} between  cohomologically Tate and mixed Tate conditions. Any hyperplane arrangement complement is mixed Tate,
 see \cite{KL}, Proposition (3.1.1).

We say that $Y$  has {\it weak polynomial count with count polynomial} $P_Y$, if there is a polynomial $P_Y\in \Z[t]$, such that for all but finitely many primes $p$, there is an integer $k_p>0$ satisfying the following: if we set $q=p^{k_p}$, then for any finite field
$\F_{q^s}$ with  $s \in \N^*$, the number of points of $Y$ over $\F_{q^s}$ is precisely $P_Y(q^s)$.

M. Kisin and G. Lehrer  have shown that if a variety $Y$ is mixed Tate then 
$Y$ has weak polynomial count (see Proposition (3.4) in \cite{KL}).

\end{rk}

\end{document}